# Neutrosophic Logic -
# A Generalization of the Intuitionistic Fuzzy Logic

Florentin Smarandache
University of New Mexico
Gallup, NM 87301, USA
smarand@unm.edu

*Abstract*

In this paper one generalizes the intuitionistic fuzzy logic (IFL) and other logics to neutrosophic logic (NL). The differences between IFL and NL (and the corresponding intuitionistic fuzzy set and neutrosophic set) are pointed out.

1. **Introduction**

The paper starts with a short paragraph on non-standard analysis because it is necessary in defining non-standard real subsets and especially the non-standard unit interval $]^-0, 1^+[$, all used by neutrosophic logic. Then a survey of the evolution of logics from Boolean to neutrosophic is presented. Afterwards the neutrosophic logic components are introduced followed by the definition of neutrosophic logic and neutrosophic logic connectors which are based on set operations. Original work consists in the definition of neutrosophic logic and neutrosophic connectors as an extension of intuitionistic fuzzy logic and the comparison between NL and other logics, especially the IFL.

2. **A Small Introduction to Non-Standard Analysis**

In 1960s Abraham Robinson has developed the *non-standard analysis*, a formalization of analysis and a branch of mathematical logic, that rigorously defines the infinitesimals. Informally, an infinitesimal is an infinitely small number. Formally, x is said to be infinitesimal if and only if for all positive integers n one has $|x| < 1/n$. Let $\varepsilon > 0$ be a such infinitesimal number. The *hyper-real number set* is an extension of the real number set, which includes classes of infinite numbers and classes of infinitesimal numbers. Let's consider the non-standard finite numbers $1^+ = 1+\varepsilon$, where "1" is its standard part and "$\varepsilon$" its non-standard part, and $^-0 = 0-\varepsilon$, where "0" is its standard part and "$\varepsilon$" its non-standard part.
Then, we call $]^-0, 1^+[$ a non-standard unit interval. Obviously, 0 and 1, and analogously non-standard numbers infinitely small but less than 0 or infinitely small but greater than 1, belong to the non-standard unit interval. Actually, by "$^-a$" one signifies a monad, i.e. a set of hyper-real numbers in non-standard analysis:
   $(^-a) = \{a-x : x \in \mathbb{R}^*, x \text{ is infinitesimal}\}$,
and similarly "$b^+$" is a monad:
   $(b^+) = \{b+x : x \in \mathbb{R}^*, x \text{ is infinitesimal}\}$.
Generally, the left and right borders of a non-standard interval $]^-a, b^+[$ are vague, imprecise, themselves being non-standard (sub)sets $(^-a)$ and $(b^+)$ as defined above.
Combining the two before mentioned definitions one gets, what we would call, a binad of "$^-c^+$":

$(^-c^+) = \{c-x : x \in \mathbb{R}^*, x \text{ is infinitesimal}\} \cup \{c+x : x \in \mathbb{R}^*, x \text{ is infinitesimal}\}$, which is a collection of open punctured neighborhoods (balls) of c.

Of course, $^-a < a$ and $b^+ > b$. No order between $^-c^+$ and c.

Addition of non-standard finite numbers with themselves or with real numbers:

$^-a + b = ^-(a + b)$

$a + b^+ = (a + b)^+$

$^-a + b^+ = ^-(a + b)^+$

$^-a + ^-b = ^-(a + b)$ (the left monads absorb themselves)

$a^+ + b^+ = (a + b)^+$ (analogously, the right monads absorb themselves).

Similarly for subtraction, multiplication, division, roots, and powers of non-standard finite numbers with themselves or with real numbers.

By extension let $\inf\ ]^-a, b^+[ = ^-a$ and $\sup\ ]^-a, b^+[ = b^+$.

## 3. A Short History

The idea of tripartition (truth, falsehood, indeterminacy) appeared in 1764 when J. H. Lambert investigated the credibility of one witness affected by the contrary testimony of another. He generalized Hooper's rule of combination of evidence (1680s), which was a Non-Bayesian approach to find a probabilistic model. Koopman in 1940s introduced the notions of lower and upper probability, followed by Good, and Dempster (1967) gave a rule of combining two arguments. Shafer (1976) extended it to the Dempster-Shafer Theory of Belief Functions by defining the Belief and Plausibility functions and using the rule of inference of Dempster for combining two evidences proceeding from two different sources. Belief function is a connection between fuzzy reasoning and probability. The Dempster-Shafer Theory of Belief Functions is a generalization of the Bayesian Probability (Bayes 1760s, Laplace 1780s); this uses the mathematical probability in a more general way, and is based on probabilistic combination of evidence in artificial intelligence.

In Lambert "there is a chance p that the witness will be faithful and accurate, a chance q that he will be mendacious, and a chance 1-p-q that he will simply be careless" [apud Shafer (1986)]. Therefore three components: accurate, mendacious, careless, which add up to 1.

Van Fraassen introduced the supervaluation semantics in his attempt to solve the sorites paradoxes, followed by Dummett (1975) and Fine (1975). They all tripartitioned, considering a vague predicate which, having border cases, is undefined for these border cases. Van Fraassen took the vague predicate 'heap' and extended it positively to those objects to which the predicate definitively applies and negatively to those objects to which it definitively doesn't apply. The remaining objects border was called penumbra. A sharp boundary between these two extensions does not exist for a soritical predicate. Inductive reasoning is no longer valid too; if S is a sorites predicate, the proposition "$\exists n(Sa_n \& \neg Sa_{n+1})$" is false. Thus, the predicate Heap (positive extension) = true, Heap (negative extension) = false, Heap (penumbra) = indeterminate.

Narinyani (1980) used the tripartition to define what he called the "indefinite set", and Atanassov (1982) continued on tripartition and gave five generalizations of the fuzzy set, studied their properties and applications to the neural networks in medicine:

  a) Intuitionistic Fuzzy Set (IFS):

  Given an universe E, an IFS A over E is a set of ordered triples <universe_element, degree_of_membership_to_A(M), degree_of_non-membership_to_A(N)> such that M+N ≤ 1 and M, N ∈ [0, 1]. When M + N = 1 one obtains the fuzzy set, and if M + N < 1 there is an indeterminacy I = 1-M-N.

  b) Intuitionistic L-Fuzzy Set (ILFS):

  Is similar to IFS, but M and N belong to a fixed lattice L.

  c) Interval-valued Intuitionistic Fuzzy Set (IVIFS):

  Is similar to IFS, but M and N are subsets of [0, 1] and sup M + sup N ≤ 1.

  d) Intuitionistic Fuzzy Set of Second Type (IFS2):

Is similar to IFS, but $M^2 + N^2 \leq 1$. M and N are inside of the upper right quarter of unit circle.
  e) Temporal IFS:
Is similar to IFS, but M and N are functions of the time-moment too.

## 4. Definition of Neutrosophic Components

Let T, I, F be standard or non-standard real subsets of $]^-0, 1^+[$,
  with   sup T = t_sup, inf T = t_inf,
         sup I = i_sup,  inf I = i_inf,
         sup F = $f_{sup}$, inf F = $f_{inf}$,
  and    $n_{sup} = t_{sup}+i_{sup}+f_{sup}$,
         $n_{inf} = t_{inf}+i_{inf}+f_{inf}$.

The sets T, I, F are not necessarily intervals, but may be any real sub-unitary subsets: discrete or continuous; single-element, finite, or (countably or uncountably) infinite; union or intersection of various subsets; etc.
They may also overlap. The real subsets could represent the relative errors in determining t, i, f (in the case when the subsets T, I, F are reduced to points).
In the next papers, T, I, F, called *neutrosophic components*, will represent the truth value, indeterminacy value, and falsehood value respectively referring to neutrosophy, neutrosophic logic, neutrosophic set, neutrosophic probability, neutrosophic statistics.

This representation is closer to the human mind reasoning. It characterizes/catches the *imprecision* of knowledge or linguistic inexactitude received by various observers (that's why T, I, F are subsets - not necessarily single-elements), *uncertainty* due to incomplete knowledge or acquisition errors or stochasticity (that's why the subset I exists), and *vagueness* due to lack of clear contours or limits (that's why T, I, F are subsets and I exists; in particular for the appurtenance to the neutrosophic sets).
One has to specify the superior (x_sup) and inferior (x_inf) limits of the subsets because in many problems arises the necessity to compute them.

## 5. Definition of Neutrosophic Logic

A logic in which each proposition is estimated to have the percentage of truth in a subset T, the percentage of indeterminacy in a subset I, and the percentage of falsity in a subset F, where T, I, F are defined above, is called *Neutrosophic Logic*.

We use a subset of truth (or indeterminacy, or falsity), instead of a number only, because in many cases we are not able to exactly determine the percentages of truth and of falsity but to approximate them: for example a proposition is between 30-40% true and between 60-70% false, even worst: between 30-40% or 45-50% true (according to various analyzers), and 60% or between 66-70% false.
The subsets are not necessary intervals, but any sets (discrete, continuous, open or closed or half-open/half-closed interval, intersections or unions of the previous sets, etc.) in accordance with the given proposition.
A subset may have one element only in special cases of this logic.

Constants: (T, I, F) truth-values, where T, I, F are standard or non-standard subsets of the non-standard interval $]^-0, 1^+[$, where $n_{inf}$ = inf T + inf I + inf F $\geq ^-0$, and $n_{sup}$ = sup T + sup I + sup F $\leq 3^+$.
Atomic formulas: a, b, c, … .
Arbitrary formulas: A, B, C, … .

Therefore, we finally generalize the intuitionistic fuzzy logic to a transcendental logic, called "neutrosophic logic": where the interval [0, 1] is exceeded, i.e. , the percentages of truth, indeterminacy, and falsity are approximated by non-standard subsets – not by single numbers, and these subsets may overlap and exceed the unit interval in the sense of the non-standard analysis; also the superior sums and inferior sum, $n_{sup}$ = sup T + sup I + sup F $\in$ ]$^-$0, 3$^+$[, may be as high as 3 or 3$^+$, while $n_{inf}$ = inf T + inf I + inf F $\in$ ]$^-$0, 3$^+$[, may be as low as 0 or $^-$0.

Let's borrow from the modal logic the notion of "world", which is a semantic device of what the world might have been like. Then, one says that the neutrosophic truth-value of a statement A, $NL_t(A) = 1^+$ if A is 'true in all possible worlds' (syntagme first used by Leibniz) and all conjunctures, that one may call "absolute truth" (in the modal logic it was named *necessary truth,* Dinulescu-Câmpina (2000) names it 'intangible absolute truth' ), whereas $NL_t(A) = 1$ if A is true in at least one world at some conjuncture, we call this "relative truth" because it is related to a 'specific' world and a specific conjuncture (in the modal logic it was named *possible truth*).
Similarly for absolute and relative falsehood and absolute and relative indeterminacy.
The neutrosophic inference ([3]), especially for plausible and paradoxist information, is still a subject of intense research today.

### 6. Differences between Neutrosophic Logic and Intuitionistic Fuzzy Logic

The differences between IFL and NL (and the corresponding intuitionistic fuzzy set and neutrosophic set) are:
   a) Neutrosophic Logic can distinguish between *absolute truth* (truth in all possible worlds, according to Leibniz) and *relative truth* (truth in at least one world), because NL(absolute truth)=1$^+$ while NL(relative truth)=1. This has application in philosophy (see the neutrosophy). That's why the unitary standard interval [0, 1] used in IFL has been extended to the unitary non-standard interval ]$^-$0, 1$^+$[ in NL.
Similar distinctions for absolute or relative falsehood, and absolute or relative indeterminacy are allowed in NL.
   b) In NL there is no restriction on T, I, F other than they are subsets of ]$^-$0, 1$^+$[, thus:
$^-$0 $\leq$ inf T + inf I + inf F $\leq$ sup T + sup I + sup F $\leq$ 3$^+$.
This non-restriction allows paraconsistent, dialetheist, and incomplete information to be characterized in NL {i.e. the sum of all three components if they are defined as points, or sum of superior limits of all three components if they are defined as subsets can be >1 (for paraconsistent information coming from different sources) or < 1 for incomplete information}, while that information can not be described in IFL because in IFL the components T (truth), I (indeterminacy), F (falsehood) are restricted either to t+i+f=1 or to $t^2 + f^2 \leq 1$, if T, I, F are all reduced to the points t, i, f respectively, or to sup T + sup I + sup F = 1 if T, I, F are subsets of [0, 1].
   c) In NL the components T, I, F can also be *non-standard* subsets included in the unitary non-standard interval
]$^-$0, 1$^+$[, not only *standard* subsets included in the unitary standard interval [0, 1] as in IFL.
   d) NL, like dialetheism, can describe paradoxes, NL(paradox) = (1, I, 1), while IFL can not describe a paradox because the sum of components should be 1 in IFL
([11],[12],[13]).
   e) NL has a better and clear name "neutrosophic" (which means the neutral part: i.e. neither true nor false), while IFL's name "intuitionistic" produces confusion with Intuitionistic Logic, which is something different.

### 7. Operations with Sets

We need to present these set operations in order to be able to introduce the neutrosophic connectors.
Let $S_1$ and $S_2$ be two (unidimensional) real standard or non-standard subsets, then one defines:

*7.1 Addition of Sets:*
$S_1 \oplus S_2 = \{x \mid x=s_1+s_2,$ where $s_1 \in S_1$ and $s_2 \in S_2\}$,

with inf $S_1 \oplus S_2$ = inf $S_1$ + inf $S_2$, sup $S_1 \oplus S_2$ = sup $S_1$ + sup $S_2$;
and, as some particular cases, we have
$\{a\} \oplus S_2 = \{x \mid x=a+s_2, \text{ where } s_2 \in S_2\}$
with inf $\{a\} \oplus S_2$ = a + inf $S_2$, sup $\{a\} \oplus S_2$ = a + sup $S_2$.

*7.2 Subtraction of Sets:*
$S_1 \ominus S_2 = \{x \mid x=s_1-s_2, \text{ where } s_1 \in S_1 \text{ and } s_2 \in S_2\}$.
For real positive subsets (most of the cases will fall in this range) one gets
inf $S_1 \ominus S_2$ = inf $S_1$ - sup $S_2$, sup $S_1 \ominus S_2$ = sup $S_1$ - inf $S_2$;
and, as some particular cases, we have
$\{a\} \ominus S_2 = \{x \mid x=a-s_2, \text{ where } s_2 \in S_2\}$,
with inf $\{a\} \ominus S_2$ = a - sup $S_2$, sup $\{a\} \ominus S_2$ = a - inf $S_2$;
also $\{1^+\} \ominus S_2 = \{x \mid x=1^+-s_2, \text{ where } s_2 \in S_2\}$,
with inf $\{1^+\} \ominus S_2 = 1^+$ - sup $S_2$, sup $\{1^+\} \ominus S_2$ = 100 - inf $S_2$.

*7.3 Multiplication of Sets:*
$S_1 \odot S_2 = \{x \mid x=s_1 \cdot s_2, \text{ where } s_1 \in S_1 \text{ and } s_2 \in S_2\}$.
For real positive subsets (most of the cases will fall in this range) one gets
inf $S_1 \odot S_2$ = inf $S_1 \cdot$ inf $S_2$, sup $S_1 \odot S_2$ = sup $S_1 \cdot$ sup $S_2$;
and, as some particular cases, we have
$\{a\} \odot S_2 = \{x \mid x=a \cdot s_2, \text{ where } s_2 \in S_2\}$,
with inf $\{a\} \odot S_2$ = a * inf $S_2$, sup $\{a\} \odot S_2$ = a · sup $S_2$;
also $\{1^+\} \odot S_2 = \{x \mid x=1 \cdot s_2, \text{ where } s_2 \in S_2\}$,
with inf $\{1^+\} \odot S_2 = 1^+ \cdot$ inf $S_2$, sup $\{1^+\} \odot S_2 = 1^+ \cdot$ sup $S_2$.

*7.4 Division of a Set by a Number:*
Let $k \in \mathbb{R}^*$, then $S_1 \oslash k = \{x \mid x=s_1/k, \text{ where } s_1 \in S_1\}$.

**8. Neutrosophic Logic Connectors**

One uses the definitions of neutrosophic probability and neutrosophic set operations.
Similarly, there are many ways to construct such connectives according to each particular problem to solve; here we present the easiest ones:

One notes the neutrosophic logic values of the propositions $A_1$ and $A_2$ by
$NL(A_1) = (T_1, I_1, F_1)$ and $NL(A_2) = (T_2, I_2, F_2)$ respectively.
For all neutrosophic logic values below: if, after calculations, one obtains numbers < 0 or > 1, one replaces them by $^-0$ or $1^+$ respectively.

*8.1. Negation:*
$NL(\neg A_1) = (\{1^+\} \ominus T_1, \{1^+\} \ominus I_1, \{1^+\} \ominus F_1)$.

*8.2. Conjunction:*
$NL(A_1 \wedge A_2) = (T_1 \odot T_2, I_1 \odot I_2, F_1 \odot F_2)$.
(And, in a similar way, generalized for n propositions.)

*8.3 Weak or inclusive disjunction:*
$NL(A_1 \vee A_2) = (T_1 \oplus T_2 \ominus T_1 \odot T_2, I_1 \oplus I_2 \ominus I_1 \odot I_2, F_1 \oplus F_2 \ominus F_1 \odot F_2)$.
(And, in a similar way, generalized for n propositions.)

*8.4 Strong or exclusive disjunction:*

$NL(A_1 \vee A_2) =$
$(T_1 \odot (\{1\} \ominus T_2) \oplus T_2 \odot (\{1\} \ominus T_1) \ominus T_1 \odot T_2 \odot (\{1\} \ominus T_1) \odot (\{1\} \ominus T_2),$
$I_1 \odot (\{1\} \ominus I_2) \oplus I_2 \odot (\{1\} \ominus I_1) \ominus I_1 \odot I_2 \odot (\{1\} \ominus I_1) \odot (\{1\} \ominus I_2),$
$F_1 \odot (\{1\} \ominus F_2) \oplus F_2 \odot (\{1\} \ominus F_1) \ominus F_1 \odot F_2 \odot (\{1\} \ominus F_1) \odot (\{1\} \ominus F_2)).$
(And, in a similar way, generalized for n propositions.)

*8.5 Material conditional (implication):*
$NL(A_1 \mapsto A_2) = (\{1^+\} \ominus T_1 \oplus T_1 \odot T_2, \{1^+\} \ominus I_1 \oplus I_1 \odot I_2, \{1^+\} \ominus F_1 \oplus F_1 \odot F_2).$

*8.6 Material biconditional (equivalence):*
$NL(A_1 \leftrightarrow A_2) = ((\{1^+\} \ominus T_1 \oplus T_1 \odot T_2) \odot (\{1^+\} \ominus T_2 \oplus T_1 \odot T_2),$
$(\{1^+\} \ominus I_1 \oplus I_1 \odot I_2) \odot (\{1^+\} \ominus I_2 \oplus I_1 \odot I_2),$
$(\{1^+\} \ominus F_1 \oplus F_1 \odot F_2) \odot (\{1^+\} \ominus F_2 \oplus F_1 \odot F_2)).$

*8.7 Sheffer's connector:*
$NL(A_1 \mid A_2) = NL(\neg A_1 \vee \neg A_2) = (\{1^+\} \ominus T_1 \odot T_2, \{1^+\} \ominus I_1 \odot I_2, \{1^+\} \ominus F_1 \odot F_2).$

*8.8 Peirce's connector:*
$NL(A_1 \downarrow A_2) = NL(\neg A_1 \wedge \neg A_2) =$
$= ((\{1^+\} \ominus T_1) \odot (\{1^+\} \ominus T_2), (\{1^+\} \ominus I_1) \odot (\{1^+\} \ominus I_2), (\{1^+\} \ominus F_1) \odot (\{1^+\} \ominus F_2)).$

## 9. Generalizations

When all neutrosophic logic set components are reduced to one element, then
$t_{sup} = t_{inf} = t$, $i_{sup} = i_{inf} = i$, $f_{sup} = f_{inf} = f$, and $n_{sup} = n_{inf} = n = t+i+f$, therefore neutrosophic logic generalizes:
- the *intuitionistic logic*, which supports incomplete theories (for $0 < n < 1$ and i=0, $0 \leq t, i, f \leq 1$);
- the *fuzzy logic* (for n = 1 and i = 0, and $0 \leq t, i, f \leq 1$);
from "CRC Concise Concise Encyclopedia of Mathematics", by Eric W. Weisstein, 1998, the fuzzy logic is "an extension of two-valued logic such that statements need not to be True or False, but may have a degree of truth between 0 and 1";
- the *intuitionistic fuzzy logic* (for n=1);
- the *Boolean logic* (for n = 1 and i = 0, with t, f either 0 or 1);
- the *multi-valued logic* (for $0 \leq t, i, f \leq 1$);
definition of <many-valued logic> from "The Cambridge Dictionary of Phylosophy", general editor Robert Audi, 1995, p. 461: "propositions may take many values beyond simple truth and falsity, values functionally determined by the values of their components"; Lukasiewicz considered three values (1, 1/2, 0). Post considered m values, etc. But they varied in between 0 and 1 only. In the neutrosophic logic a proposition may take values even greater than 1 (in percentage greater than 100%) or less than 0.
- the *paraconsistent logic*, which support conflicting information (for n > 1 and i = 0, with both t, f < 1);
the *dialetheism*, which says that some contradictions are true (for t = f = 1 and i = 0; some paradoxes can be denoted this way too);
the *faillibilism*, which says that uncertainty belongs to every proposition (for i > 0);
Compared with all other logics, the neutrosophic logic and intuitionistic fuzzy logic introduce a percentage of "indeterminacy" - due to unexpected parameters hidden in some propositions, or unknowness, but neutrosophic logic let each component t, i, f be even boiling *over 1* (overflooded), i.e. be $1^+$, or freezing *under 0* (underdried), i.e. be $^-0$ in order to be able to make distinction between relative truth and absolute truth, and between relative falsity and absolute falsity in philosophy.